\documentclass[reqno]{amsart}
\theoremstyle{plain}
\newtheorem{thm}{\bf Theorem}[section]

\theoremstyle{remark}
\newtheorem{defn}[thm]{\bf Def{}inition}
\newtheorem{rem}[thm]{\bf Remark}
\newtheorem{exa}[thm]{\bf Example}
\numberwithin{equation}{section}
\usepackage{enumerate}

\allowdisplaybreaks
\begin{document}
\baselineskip12pt
\title[Wave Packet Frames ]{ Wave packet frames  generated by hyponormal operators on $L^2(\mathbb{R})$}

\author[ Vashisht]{Lalit  Kumar  Vashisht\\
\textbf{PRINCIPAL INVESTIGATOR}}

\address{L. K. Vashisht, Department of Mathematics,
University of Delhi, Delhi-110007, India}
\email{ lalitkvashisht@gmail.com}

\begin{abstract}\baselineskip8pt
 In this paper we study frame-like properties of a  wave packet system   by using  hyponormal   operators on $L^2(\mathbb{R})$. We present necessary and sufficient conditions in terms of relative hyponormality of operators for a system to be a wave packet frame in $L^2(\mathbb{R})$.  A characterization of hyponormal operators by using  tight wave packet frames is proved. This is different from a method proved by Djordjevi$\acute{c}$ by using the Moore-Penrose inverse of a bounded linear operator with a closed range.  We extends some results by Kaushik, Singh and Virender to wave packet  frames generated by hyponormal operators
 .
\end{abstract}

  \keywords{ Wave packet system, analysis operator, frame operator, hyponormal operator,  Hilbert space frame.\\
Lalit  was    supported by R $\&$ D Doctoral Research Programme, University of Delhi,
Delhi-110007, India. Grant No. : RC/2014/6820.}

\maketitle \thispagestyle{empty} \baselineskip8pt

\maketitle
\section{Introduction}\baselineskip12pt
Frames in Hilbert spaces are a redundant system  of vectors which provides a series representation for
each vector in the space. Duffin and  Schaffer \cite{DS} in 1952, introduced frames for Hilbert spaces, in the context of nonharmonic
Fourier series. Frames were revived by Daubechies, Grossmann and Meyer in \cite{DGM}. For applications of frames in various directions, see \cite{CK, OC}

Feichtinger and Werther  \cite{FW}  introduced a family of analysis and synthesis systems with frame-like properties  for closed subspaces of a separable Hilbert space $\mathcal{H}$ and call it an \emph{atomic system} (or \emph{local atoms}). The motivation for the atomic system is based on examples arising in sampling theory. One of the important properties of the atomic system is that it can generate a proper subspace even though they do not belong to them.
 \begin{defn}\cite{FW}
 Let $\mathcal{H}$ be a Hilbert space and let $\mathcal{H}_0$ be a closed subspace of $\mathcal{H}$. A sequence $\{f_k\} \subset
\mathcal{H}$ is called a \emph{family of local atoms} (or \emph{atomic system}) for
$\mathcal{H}_0$, if
\begin{enumerate}[$(i)$]
 \item there exists  a real number  $B > 0$ such that  $\|\{\langle f, f_k\rangle\}\|^2_{\ell^2} \leq B \|f\|^2$
for  all $f \in \mathcal{H}$,
\item there exists a sequence of linear functionals $\{c_k\}$ and  a real number  $C > 0$ such that
\begin{align*}
 \|\{c_k(f)\}\|^2_{\ell^2} \leq C \|f\|^2 \ \text{for  all} \  f \in \mathcal{H}_0
\intertext{and}
f = \sum_{k=1}^{\infty}c_k(f) f_k \ \text{for  all} \ f \in \mathcal{H}_0.
\end{align*}
\end{enumerate}

 \end{defn}

 G$\check{a}$vruta in \cite{LG}  introduced and studied $K$-frames
in Hilbert spaces to study atomic systems with respect to a bounded
linear operator $K$ on Hilbert spaces.
\begin{defn}\cite{LG}
Let $\mathcal{H}$ be a Hilbert space and let $K$ be a bounded linear
operator on $\mathcal{H}$. A sequence $\{f_k\} \subset
\mathcal{H}$ is called a \emph{$K$-frame} for $\mathcal{H}$, if there exist
 constants $A, B
> 0$ such that
\begin{align}
A\|K^*f\|^2\leq \sum_{k=1}^{\infty} |\langle f, f_k\rangle|^2 \leq
B\|f\|^2 \ \text{for all} \ f \in \mathcal{H}.
\end{align}
\end{defn}
  The lower inequality in $(1.1)$ is controlled by a bounded linear operator on $\mathcal{H}$. It is observed in \cite{LG}   that
$K$-frames are more general than standard frames in the sense that
the lower frame bound only holds for the elements in the range of
$K$, where $K$ is a  bounded linear operator on the underlying
Hilbert space. G$\check{a}$vruta in \cite{LG}  characterize $K$-frames in
Hilbert spaces by using bounded linear operators.

It would be interesting to control both lower and upper frame condition in $(1.1)$ by bounded linear operators on  $\mathcal{H}$. In this direction, we study  frame-like properties of an irregular wave packet system in $L^2(\mathbb{R})$, where both lower and upper frame conditions are controlled by bounded linear operators on $L^2(\mathbb{R})$ (see Definition 3.1).  The wave packet system is a family of  functions generated by combined  action of dilation, translation and modulation operators on  $L^2(\mathbb{R})$. More precisely, we consider a system of the form
\begin{align}
\{D_{a_j}T_{bk}E_{c_m}\psi\}_{j, k, m\in \mathbb{Z}},
\end{align}
where  $\psi\in L^2(\mathbb{R})$, $\{a_j\}_{j\in \mathbb{Z}}\subset
\mathbb{R}^{+}$, $b \ne 0$ and  $\{c_m\}_{m\in \mathbb{Z}} \subset \mathbb{R}$  and  call it \emph{irregular Weyl-Heisenberg wave packet system} (or simply\emph{ wave packet system}) in  $L^2(\mathbb{R})$.  A frame for  $L^2(\mathbb{R})$ of the form  $\{D_{a_j}T_{bk}E_{c_m}\psi\}_{j, k, m\in \mathbb{Z}}$ is called an
\emph{irregular wave packet frame} (or\emph{ wave packet frame}).  The  wave packet system was introduced  by Cordoba
and Fefferman \cite{CF} by applying certain collections of dilations,
modulations and translations to the Gaussian function in the study
of some classes of singular integral operators. Later, Labate et al.
\cite{LWW} adopted the same expression to describe, more generally,  any collection of
functions which are obtained by applying the same operations to a
finite family of functions in $L^2(\mathbb{R}^d)$. More precisely, Gabor systems, wavelet
systems and the Fourier transform of wavelet systems are special
cases of wave packet systems. Lacey and Thiele  \cite{LT1, LT2} gave applications of wave packet systems in boundedness of the  Hilbert transforms. The  wave packet systems have been studied by several authors, see \cite{CKS, GL, HLW, HLWW, SLV, SV}.

\subsection{Outline:}
This paper is organized as follows: In Section 2, we give basic
definitions and results which will be used throughout the paper.
 Section 3 is devoted to the study of frame-like properties of irregular Weyl-Heisenberg wave packet  systems. We introduce  \emph{$\Theta$-irregular Weyl-Heisenberg wave packet frame} (in short, \emph{$\Theta$-$I W H$ wave packet frame}) for $L^2(\mathbb{R})$, where $\Theta$ is a bounded linear operator on $L^2(\mathbb{R})$ (see  Definition 3.1). This type of  wave packet frame can control both lower and upper frame conditions by bounded linear operators on $L^2(\mathbb{R})$ . The $\Theta$-$IWH$ wave packet frame (in the context of standard Hilbert frame) for a  Hilbert space is a $K$-frame, but converse is not true (see Example 3.2). Furthermore, the  $\Theta$-$I W H$ wave packet frame control both lower and upper frame conditions by bounded linear operators.  Necessary and sufficient conditions for a certain system to be a $\Theta$-$IWH$ wave packet frames  for $L^2(\mathbb{R})$  by using hyponormality of  operators  on $L^2(\mathbb{R})$ have  been obtained. A characterization of hyponormal operator  in terms of a special  type of tight wave packet frames for $L^2(\mathbb{R})$ is given. This is different from a method proved by Djordjevi$\acute{c}$ in \cite{DJ} by using the Moore-Penrose inverse of a bounded linear operator with a closed range (see Theorem 3.7).   The  linear combinations of frames or redundant building blocks are important in applied mathematics, we  discuss linear combinations of $\Theta$-$IWH$ wave packet frames for $L^2(\mathbb{R})$ in Section 4.

\section{Preliminaries}
In this section, we recall basic notations and definitions to  make the paper self-contained.
Let $\mathcal{H}$ be a separable real (or complex)  Hilbert space  with inner product $\langle ., . \rangle$ linear in the first entry.  A countable sequence $\{f_k\} \subset \mathcal{H}$ is
called a \emph{frame} \break (or \emph{Hilbert frame}) for $\mathcal{H}$, if there exist numbers $0 <  a_o \leq b_o < \infty$ such that
\begin{align}
a_o \|f\|^2\leq  \sum_{k=1}^{\infty} |\langle f, f_k\rangle|^2 \leq b_o \|f\|^2
\ \text{for  all} \  f \in \mathcal{H}.
\end{align}
The numbers $a_o$ and $b_o$ are called \emph{lower} and \emph{upper frame bounds}, respectively. They are not unique. If  it is possible to choose $a_o = b_o$, then the frame $\{f_k\}$ is called \emph{Parseval frame }(or \emph{tight frame}).\\
The scalars
\begin{align*}
& \mathrm{\gamma}_0= \inf \{b_o > 0: b_o  \text{ satisfies (2.1)}\}\\
&\mathrm{\delta}_0= \sup \{a_o > 0: a_o  \text{ satisfies (2.1)}\}
\end{align*}
are called the \emph{optimal bounds} or \emph{best bounds} of the  frame.

Associated with a frame $\{f_k\}$  for $\mathcal{H}$, there are three bounded linear operators:
\begin{align*}
& \text{\emph{synthesis operator}} \quad  V:\ell^2\rightarrow \mathcal{H}, \quad  \ V(\{c_k\})=\sum_{k=1}^{\infty}c_kf_k,  \ \{c_k\} \in \ell^2,\\
& \text{\emph{analysis operator}} \quad  V^{*}:\mathcal{H}\rightarrow
\ell^2,  \quad   V^{*}(f)=\{\langle f, f_k\rangle\},\  \ f \in
\mathcal{H},\\
& \text{\emph{frame operator}} \quad S=V V^{*}:\mathcal{H}\rightarrow
\mathcal{H}, \quad   S(f)=\sum_{k=1}^{\infty}\langle f, f_k\rangle f_k,  \ \
f \in \mathcal{H}.
\end{align*}

 The frame operator $S$ is a positive,
self-adjoint and invertible operator on $\mathcal{H}$. This gives the \emph{reconstruction formula} for all $f \in \mathcal{H}$,
\begin{align*}
&f = SS^{-1}f =\sum_{k=1}^{\infty}\langle S^{-1}f, f_k \rangle f_k
 \quad  \left(\ =\sum_{k=1}^{\infty}\langle f,S^{-1}f_k\rangle f_k \ \right).
\end{align*}
The scalars $\{\langle S^{-1}f, f_k \rangle\}$ are called \emph{frame coefficients} of the vector $f \in \mathcal{H}$.  The representation of
$f$ in the reconstruction formula need not be unique. This reflects one of the important properties of frames in applied mathematics.

Let $a, b\in \mathbb{R}$ and $c \in \mathbb{R} \setminus \{0\}$. We consider  operators  $T_a, \ E_b, \ D_c : L^2(\mathbb{R})  \rightarrow~L^2(\mathbb{R})$ given by
\begin{align*}
& \text{Translation by} \ a \leftrightarrow  \ T_af(t)=f(t-a),\\
& \text{Modulation by} \ b \leftrightarrow  \ E_bf(t)=e^{2\pi ib t}f(t), \\
& \text{Dilation by} \  c \leftrightarrow D_{c}f(t) = |c|^{\frac{1}{2}} f(ct).
\end{align*}
 A bounded linear operator $T$  defined on
$\mathbb{H}$ is said to be \emph{positive}, if $\langle Tf, f
\rangle \geq 0$ for all $f \in \mathbb{H}$. In symbol we write $T
\geq 0$. If $T_1, T_2$  are bounded linear operator on $\mathbb{H}$ such that $T_1 - T_2
\geq 0$, then we write $T_1\geq T_2$.  A bounded linear operator $T: \mathbb{H} \rightarrow \mathbb{H}$
 is said to be \emph{hyponormal}, if $T^*T- TT^* \geq 0$, or equivalently if $\|T^*f\| \leq \|Tf\|$ for all $f \in \mathbb{H}$. The \emph{characteristic function} of any set $E$ is denoted by $\chi_E$. By $\mathcal{R}(T)$ we denote the range of a bounded linear operator $T$ from a normed space $X$ into a normed space $Y$.

\begin{thm} \cite{D}\label{thm 2.4}
Let $\mathbb{H}, \mathbb{H}_1, \mathbb{H}_2$ be  Hilbert spaces. Assume that $T_1: \mathbb{H}_1 \rightarrow \mathbb{H}$ and  $T_2:\mathbb{H}_2 \rightarrow \mathbb{H}$ be
bounded linear operators.
The following statement are equivalent:
\begin{enumerate}[$(i)$]
\item $\mathcal{R}(T_1)\subset  \mathcal{R}(T_2)$.
\item $T_1T_1^{*}\leq \lambda ^2 T_2T_2^{*}$ for some $\lambda\geq0$.
\item There exists a bounded linear operator $S:\mathbb{H}_1 \rightarrow \mathbb{H}_2$ such
that $T_1=T_2S .$
\end{enumerate}
\end{thm}
\section{Wave Packet Frames in  $L^2(\mathbb{R})$}
\begin{defn}
Let $\psi\in L^2(\mathbb{R})$, $\{a_j\}_{j\in \mathbb{Z}}\subset
\mathbb{R}^+$,  $\{c_m\}_{m\in \mathbb{Z}} \subset
\mathbb{R}$ and $b\neq0$ and let $\Theta$ be a bounded linear operator on $L^2(\mathbb{R})$. A system $\{D_{a_j}T_{bk}E_{c_m}\psi\}_{j,k,m \in
\mathbb{Z}} $ is called a \emph{$\Theta$-irregular Weyl-Heisenberg wave packet frame} (in short, \emph{$\Theta$-$I W H$ wave packet frame}) for $L^2(\mathbb{R})$,
 if there exist   constants
$ 0~<~\mathrm{\alpha_0} \leq \mathrm{\beta_0} < \infty$ such that
\begin{align}
\mathrm{\alpha_0}\|\Theta^{*}f\|^2\leq \sum_{j,k,m \in
\mathbb{Z}}|\langle
 f,D_{a_j}T_{bk}E_{c_m}\psi\rangle|^2 \leq \mathrm{\beta_0}\|\Theta f\|^2 \  \text{for   all} \ f\in L^2(\mathbb{R}).
 \end{align}
\end{defn}

The scalars $\alpha_0$ and $\beta_0$ are called \textit{lower}
 and \emph{upper  bounds} of the $\Theta$-$I W H$ wave packet frame $\{D_{a_j}T_{bk}E_{c_m}\psi\}_{j,k,m \in \mathbb{Z}} $,
respectively. If upper inequality  in $(3.1)$ is satisfied, then
$\{D_{a_j}T_{bk}E_{c_m}\psi\}_{j,k,m \in \mathbb{Z}} $ is called a
\emph{ Bessel sequence} in
$L^2(\mathbb{R})$ with Bessel bound $\mathrm{\beta_0}$.
If $\Theta$ is the identity operator on $L^2(\mathbb{R})$, then $\Theta$-$I W H$ wave packet frame for $L^2(\mathbb{R})$ is the standard $I W H$ wave packet frame for $L^2(\mathbb{R})$.

If a countable sequence $\{f_k\}$ in a Hilbert space $\mathcal{H}$  satisfies the inequality (3.1), i.e., if
\begin{align*}
\mathrm{\alpha_0}\|\Theta^{*}f\|^2\leq \sum_{k=1}^{\infty}|\langle f, f_k\rangle|^2 \leq \mathrm{\beta_0}\|\Theta f\|^2 \  \text{for   all} \ f \in \mathcal{H},
 \end{align*}
then we say that $\{f_k\}$ is a \emph{$\Theta$-Hilbert frame} for $\mathcal{H}$.

\subsection{Examples and comments:}
 Every $\Theta$-Hilbert frame for $\mathcal{H}$ is a $K$-frame for $\mathcal{H}$, but not conversely. More precisely, if $\{f_k\}$ is a $\Theta$-Hilbert frame for $\mathcal{H}$  with frame  bounds $\alpha_0$ and $\beta_0$. Then, $\{f_k\}$ is a $K$-frame for $\mathcal{H}$ with frame bounds $\alpha_0$ and $\beta_0 \ \|\Theta \|^2$. The following example shows that a $K$-frame for $\mathcal{H}$ need not be a $\Theta$-Hilbert frame for $\mathcal{H}$.

\begin{exa}
Let $\{\chi_k\}$ be the canonical orthonormal basis for the discrete signal space $\mathcal{H} = L^2(\Omega, \mu)$ (where $\Omega = \mathbb{N}$ and $\mu$ is the counting measure) and let $\Theta$ be the backward shift operator on $\mathcal{H}$ given by
\begin{align*}
\Theta( \{\xi_1, \xi_2, \xi_3,....... \}) = \{\xi_2,  \xi_3,....... \}, \ \{\xi_j\} \in \mathcal{H}.
\end{align*}
Then, its conjugate $\Theta^*$ is the forward shift operator on $\mathcal{H}$ which is given by
\begin{align*}
\Theta^*( \{\xi_1, \xi_2, \xi_3,....... \}) = \{0, \xi_1, \xi_2,  \xi_3,....... \}, \ \{\xi_j\} \in \mathcal{H}.
\end{align*}
Choose $f_k = \chi_k$ for all $k \in \mathbb{N}$.\\
We compute
\begin{align*}
\|\Theta^* f \|^2 = \| f \|^2 = \sum_{j=1}^{\infty} |\langle f,  f_k \rangle|^2 \ \text{for all} \ f =  \{\xi_j\} \in \mathcal{H}.
\end{align*}
Hence $\{f_k\}$ is a $K$-frame (with a choice $K = \Theta)$ for $\mathcal{H}$ with frame bounds $A = B = 1$. But $\{f_k\}$ is not a $\Theta$-Hilbert frame for $\mathcal{H}$.
Indeed, let $a_o$ and $b_o$ be positive numbers such that
\begin{align}
a_o\|\Theta^*f\|^2\leq \sum_{k=1}^{\infty} |\langle f, f_k\rangle|^2 \leq b_o\|\Theta f\|^2 \ \text{for all} \ f \in \mathcal{H}.
\end{align}
Then, for  $f_o=\chi_1 \in \mathcal{H}$, we obtain $\Theta f_o=0$. Therefore, by using upper inequality  in $(3.2)$, we have  $f_o = 0$, a contradiction.
\end{exa}

\begin{rem}
A $\Theta$-Hilbert frame for $\mathcal{H}$  ($\Theta\neq I$, the identity operator on $\mathcal{H}$) need not be a standard Hilbert frame for $\mathcal{H}$  and vice-versa. Indeed, let $\mathcal{H}$ be the  discrete signal space given in Example 3.2  with canonical orthonormal basis $\{\chi_k\}$.\\
 Choose $f_k = \chi_k + \chi_{k+1},  k \in \mathbb{N}$.\\
  Define $\Theta: \mathcal{H} \rightarrow \mathcal{H}$ by
\begin{align*}
\Theta(f =  \{\xi_1, \xi_2, \xi_3,....... \}) = \{\xi_1,  \xi_1 + \xi_2, \xi_2 + \xi_3,....... \}, \ f = \{\xi_j\} \in \mathcal{H}.
\end{align*}
Then, $\Theta$ is a bounded linear operator on $\mathcal{H}$ and its conjugate operator $\Theta^*$ is given by
\begin{align*}
\Theta^*( \{\xi_1, \xi_2, \xi_3,....... \}) = \{  \xi_1 + \xi_2, \xi_2 + \xi_3,....... \}, \ \{\xi_j\} \in \mathcal{H}.
\end{align*}
 One can  verify that there exists a $\gamma \in (0, 1)$ such that
\begin{align*}
\gamma \ \|\Theta^* f\|^2 \leq  \sum_{j=1}^{\infty} |\langle f,  f_k \rangle|^2 \leq \|\Theta f\|^2 \ \text{for all} \ f  \in \mathcal{H}.
\end{align*}
Hence $\mathcal{F}  \equiv \{f_k\}$ is a $\Theta$-Hilbert frame for $\mathcal{H}$. But $\mathcal{F}$ is not a standard Hilbert frame for $\mathcal{H}$  (see Example $5.4.6$ in \cite{OC}, \ p. \  98).\\
To show that a standard Hilbert frame for $\mathcal{H}$ need not be $\Theta$-Hilbert frame for $\mathcal{H}$. Choose $g_k = \chi_k, k \in \mathbb{N}$ and let $\Theta$ be the backward shift operator on $\mathcal{H}$. Then, $\mathcal{G} = \{g_k\}$ is a Hilbert frame for $\mathcal{H}$, but not a $\Theta$-Hilbert frame for $\mathcal{H}$.
\end{rem}

Regarding the existence of  $\Theta$-$IWH$ wave packet frames for  $L^2(\mathbb{R})$, we have following examples.
\begin{exa}
Let $a > 1$ and $b > 0$ and $c_m =0$ for all $m \in \mathbb{Z}$. Choose $a_j =  a^{j}$  for all $j \in \mathbb{Z}$. Then, there exist $\psi\in L^2(\mathbb{R})$ such that
 $\hat{\psi}=\chi_E$,  where $E$ is a compact subset of $\mathbb{R}$. Therefore,
 \begin{align*}
 \{D_{a_j}T_{bk}E_{c_m}\psi\}_{j, k, m\in{\mathbb{Z}}}=\{D_{a^j}T_{bk}\psi\}_{j, k \in{\mathbb{Z}}}
 \end{align*}
 is an orthonormal basis for  $L^2({\mathbb{R})}$ (see  Theorem 12.3 in \cite{H} \ \ p. 357), hence a tight  $IWH$ wave packet  frame for $L^2(\mathbb{R})$ .\\
Let $\beta\in\mathbb{R}$ be arbitrary, but fixed.\\
  Choose $\Theta=E_{\beta}$ (the modulation operator on $L^2(\mathbb{R})$) and $d_m=c_m+\beta \ (m\in\mathbb{Z})$.\\
We compute
\begin{align*}
\| \Theta^*f \|^2 =  \alpha_o\|E_{\beta}^*f\|^2  &= \sum_{j,k,m \in \mathbb{Z}}|\langle E_{\beta}^*f,D_{a_j}T_{bk}E_{c_m}\psi\rangle|^2\\
&=\sum_{j,k,m \in \mathbb{Z}}|\langle f,D_{a_j}T_{bk}E_{d_m}\psi\rangle|^2\\
& = \| \Theta f \|^2, \ \text{for all} \ f \in L^2(\mathbb{R}).
\end{align*}
Hence $\{D_{a_j}T_{bk}E_{d_m}\psi\}_{j,k,m \in \mathbb{Z}}$ is a $\Theta$-$IWH$ wave packet frame for
$L^2(\mathbb{R})$.
\end{exa}

\begin{exa}
Let $\Theta: L^2(\mathbb{R}) \rightarrow L^2(\mathbb{R})$ be the multiplication operator given by
\begin{align*}
\Theta(f) = f . \chi_{[0,1]}, \ f \in L^2(\mathbb{R}).
\end{align*}
Then, $\Theta$ is a bounded linear operator on $L^2(\mathbb{R})$.\\
 Choose $b=1, \ a_j = 1, \  c_m =0$ \ $(j,m\in \mathbb{Z})$ and   $\psi = \chi_{[0,1]}$.\\
 Then
  \begin{align*}
  \{D_{a_j}T_{bk}E_{c_m}\psi\}_{j,k,m \in \mathbb{Z}} = \{ T_k\psi\}_{k \in \mathbb{Z}} =  \{\chi_{[k, k+1]} \}_{k \in \mathbb{Z}}.
  \end{align*}
  The system $\{D_{a_j}T_{bk}E_{c_m}\psi\}_{j,k,m \in \mathbb{Z}}$ is not a $\Theta$-$IWH$ wave packet frame for $L^2(\mathbb{R})$. Indeed, let $B$ be an upper $\Theta$-$IWH$ wave packet frame bound for $\{D_{a_j}T_{bk}E_{c_m}\psi\}_{j,k,m \in \mathbb{Z}}$.
Let $h \in L^2(\mathbb{R})$ be a function given by
 \begin{align*}
h(x)&=\begin{cases}
 \chi_{[0,1]},          & x \in [0, 1]\\
\sqrt{B} \ \ \chi_{[2,3]},  & x \in [2, 3]\\
0& \text{otherwise}.
\end{cases}
\end{align*}
 We compute
 \begin{align*}
 \sum_{j,k,m \in \mathbb{Z}}|\langle h, D_{a_j}T_{bk}E_{c_m}\psi\rangle|^2 & = \sum_{k \in \mathbb{Z}}|\langle h, \chi_{[k, k+1]}\rangle|^2\\
 & =  |\langle h, \chi_{[0, 1]}\rangle|^2 + |\langle h, \chi_{[2, 3]}\rangle|^2\\
 & = 1 + B.
\end{align*}
 On the other hand, $\|\Theta h\|^2 = \| h . \chi_{[0, 1]}\|^2 = 1$.\\
  Therefore, $\sum_{j,k,m \in \mathbb{Z}} |\langle h, D_{a_j}T_{bk}E_{c_m}\psi\rangle|^2 =  1 + B > B \|\Theta h\|^2$. Hence $B$ is not an upper  $\Theta$-$IWH$ wave packet frame bound for $\{D_{a_j}T_{bk}E_{c_m}\psi\}_{j,k,m \in \mathbb{Z}}$, a contradiction.
\end{exa}

\subsection{Operators associated with $\Theta$-$I W H$ wave packet frames}
 Suppose that $\mathcal{F} \equiv \{D_{a_j}T_{bk}E_{c_m}\psi\}_{j,k,m \in
\mathbb{Z}} $ is a $\Theta$-$I W H$ wave packet frame for $L^2(\mathbb{R})$. The
operator $T: \ell^2(\mathbb{Z}^3) \rightarrow  L^2(\mathbb{R})$ given
by
\begin{align*}
T\{c_{jkm}\}_{j,k,m \in \mathbb{Z}}=\sum_{j,k,m \in
\mathbb{Z}}c_{jkm}D_{a_j}T_{bk}E_{c_m}\psi,
 \end{align*}
is called the\textit{ pre-frame operator} or \textit{synthesis
operator} associated with $\mathcal{F}$ and the adjoint operator
$T^{*}:L^2(\mathbb{R})\rightarrow \ell^2(\mathbb{Z}^3)$ is given by
\begin{align*}
T^{*}\textit{f}=\{\langle f,D_{a_j}T_{bk}E_{c_m}\psi\rangle\}_{j,k,m \in
\mathbb{Z}}\end{align*}
 is called the
\textit{analysis operator} associated with $\mathcal{F}$. Composing
$T$  and $T^{*}$, we obtain the \textit{frame operator}
$\mathcal{S}:L^2(\mathbb{R})\rightarrow L^2(\mathbb{R})$ given
by
\begin{align}\label{eq 3.3}
\mathcal{S}f=TT^{*}f=\sum_{j,k,m \in \mathbb{Z}}\langle
 f,D_{a_j}T_{bk}E_{c_m}\psi\rangle D_{a_j}T_{bk}E_{c_m}\psi.
\end{align}
Since $\mathcal{F}$ is a  $\Theta$-$IWH$ wave packet Bessel sequence in
$L^2(\mathbb{R})$, the series defining $\mathcal{S}$ converges
unconditionally for all $f\in L^2(\mathbb{R}).$  Notice that, in
general, frame operator of the $\Theta$- $IWH$ wave packet frame $\mathcal{F}$ is not invertible on
$L^2(\mathbb{R})$, but it is invertible on a subspace
$\mathcal{R}(\Theta)\subset L^2(\mathbb{R})$. In fact, if $\mathcal{R}(\Theta)$
is closed , then there exist a pseudoinverse $\Theta^{\dagger}$ of
$\Theta$ such that $\Theta\Theta^{\dagger}f=f$ for all $f \in
\mathcal{R}(\Theta)$, i.e.,
$\Theta\Theta^{\dagger}|_{\mathcal{R}(\Theta)}=I_{\mathcal{R}(\Theta)}$,
so we have
$\left(\Theta^{\dagger}|_{\mathcal{R}(\Theta)}\right)^*\Theta^{*}=I_{\mathcal{R}(\Theta)}^{*}$.
Hence for any $f\in \mathcal{R}(\Theta)$, we obtain
\begin{align*}
\|f\|=\left\|\left(\Theta^{\dagger}|_{\mathcal{R}(\Theta)}\right)^*\Theta^{*}f\right\|\leq\|\Theta^{\dagger}\|\|\Theta^{*}f\|.
\end{align*}
Therefore, by using  \eqref{eq 3.3}, we can write
\begin{align*}
\langle \mathcal{S}f,f\rangle\geq A\|\Theta^{*}f\|^2\geq
A\|\Theta^{\dagger}\|^{-2}\|f\|^2 \ \text{ for  all } \ f \in
\mathcal{R}(\Theta).
\end{align*}
That is
\begin{align*}
A\|\Theta^{\dagger}\|^{-2}\|f\|^2\leq\|Sf\|^2\leq B \|f\|^2 \ \text{for  all} \ f \in \mathcal{R}(\Theta).
\end{align*}
Thus, the operator  $\mathcal{S}:\mathcal{R}(\Theta)\rightarrow
\mathcal{S}(\mathcal{R}(\Theta))$ is a homeomorphism. Furthermore,
we have
\begin{align*}
B^{-1}\|f\|\leq\|\mathcal{S}^{-1}f\|\leq
A^{-1}\|\Theta^{\dagger}\|^2\|f\| \  \text{for  all }   f \in
\mathcal{S}(\mathcal{R}(\Theta)).
\end{align*}
\vspace{.3cm}

Next, we characterizes  a system $\{D_{a_j}T_{bk}E_{c_m}\psi\}_{j,k,m \in \mathbb{Z}} \subset L^2(\mathbb{R})$ as $\Theta$-$IWH$ wave packet frame.
 Let $T_1: \mathbb{H} \rightarrow \mathbb{H}$ and  $T_2:\mathbb{H}_1 \rightarrow \mathbb{H}$ be bounded linear operators, where $\mathbb{H}, \mathbb{H}_1$  are Hilbert spaces.   We say that  the pair $(T_1, T_2)$ is \emph{relatively hyponormal}, if
 \begin{align*}
 \lambda T_1^* T_1 \geq T_2 \ T_2^* \ \text{for some} \ \lambda > 0.
 \end{align*}
In this case we say that $T_1$ and   $T_2$  are \emph{relatively hyponormal}. Aldroubi in \cite{A} characterized operators on  a Hilbert space $\mathcal{H}$, which can generate Hilbert frames  (as images of given frames) for $\mathcal{H}$. Actually,  Aldroubi   considered operators which are relative hyponormal with the identity operator on  $\mathcal{H}$. The following theorem characterizes  a certain system as a  $\Theta$-$IWH$ wave packet frame for $L^2(\mathbb{R})$ in terms of the relative hyponormality of operators.
\begin{thm}\label{thm 3.11}
 Let $\psi\in L^2(\mathbb{R})$, $\{a_j\}_{j\in \mathbb{Z}}\subset
\mathbb{R}^+$,  $\{c_m\}_{m\in \mathbb{Z}} \subset
\mathbb{R}$ and $b\neq0$ and let $\Theta$ be a bounded linear operator on $L^2(\mathbb{R})$. Then,
$\{D_{a_j}T_{bk}E_{c_m}\psi\}_{j,k,m \in \mathbb{Z}} $ is a $\Theta$-$IWH$ wave packet frame for $L^2(\mathbb{R})$ if and only if there
exist a  bounded linear operator  $\Xi:\ell^2(\mathbb{Z}^3) \rightarrow
L^2(\mathbb{R})$ such that
\begin{enumerate}[$(i)$]
 \item the pair $(\Theta, \Xi)$  is relative hyponormal, i.e.,  $\lambda \Theta^* \Theta \geq \Xi \ \Xi^* \ \  \text{for some}\ \lambda~>~0$,
\item $\Xi (e_{j,k,m})=D_{a_j}T_{bk}E_{c_m}\psi$ \ $(j,k,m \in \mathbb{Z})$ \ and  $\mathcal{R}(\Theta)\subset \mathcal{R}(\Xi)$,
\end{enumerate}
where $\{e_{j,k,m}\}_{j,k,m \in \mathbb{Z}}$ is an orthonormal basis for $\ell^2(\mathbb{Z}^3) $.
\end{thm}

\proof
 Suppose first  that $\{D_{a_j}T_{bk}E_{c_m}\psi\}_{j,k,m \in \mathbb{Z}}$ is a $\Theta$-$IWH$ wave packet frame
for $L^2(\mathbb{R})$. Then, we can find positive constants
$\mathrm{a}_0, \mathrm{b}_0$ such that
\begin{align}\label{eq 3.7}
\mathrm{a}_0\|\Theta^{*}f\|^2\leq \sum_{j,k,m \in \mathbb{Z}}|\langle
 f,D_{a_j}T_{bk}E_{c_m}\psi\rangle|^2 \leq \mathrm{b}_0\|\Theta f\|^2 \  \text{for  all} \  f \in L^2(\mathbb{R}).
\end{align}
Define $ \mathcal{W}:L^2(\mathbb{R}) \rightarrow \ell^2(\mathbb{Z}^3)$ by
\begin{align*}
\mathcal{W}(f)=\sum_{j,k,m \in \mathbb{Z}}\langle
 f,D_{a_j}T_{bk}E_{c_m}\psi \rangle e_{j,k,m}.
 \end{align*}
Clearly, $\mathcal{W}$ is a well defined bounded linear
operator on $L^2(\mathbb{R})$.\\
\eject
We compute
\begin{align*}
  \langle \mathcal{W}^{*}e_{j,k,m},h\rangle &= \langle e_{j,k,m},\mathcal{W}h\rangle\\
 &=\left \langle e_{j,k,m},
 \sum_{j,k,m \in \mathbb{Z}}\langle h,D_{a_j}T_{bk}E_{c_m}\psi\rangle
 e_{j,k,m}\right \rangle\\
 &=\sum_{j,k,m \in \mathbb{Z}}\overline{\langle h, D_{a_j}T_{bk}E_{c_m}\psi\rangle}\langle
 e_{j,k,m},e_{j,k,m}\rangle\\
 &=\overline{\langle h,D_{a_j}T_{bk}E_{c_m}\psi\rangle}\\
 &=\langle D_{a_j}T_{bk}E_{c_m}\psi, h\rangle  \ \text{for  all} \ h \in L^2(\mathbb{R}).
\end{align*}
This gives
\begin{align}\label{eq 3.8}
\mathcal{W}^{*}e_{j,k,m}=D_{a_j}T_{bk}E_{c_m}\psi \ \ (j,k,m \in \mathbb{Z}).
 \end{align}
By using \eqref{eq 3.8} and lower frame inequality in \eqref{eq 3.7}, we obtain
\begin{align*}
\mathrm{a}_0\|\Theta^{*}f\|^2\leq\sum _{j,k,m \in \mathbb{Z}}|\langle
f,\mathcal{W}^{*}e_{j,k,m}\rangle|^2=\|\mathcal{W}f\|^2  \text{ for  all} \ f  \in L^2(\mathbb{R}).
\end{align*}
This gives $\mathrm{a}_0\Theta\Theta^{*}\leq \mathcal{W}^{*}\mathcal{W}$.\\
Choose $\Xi =\mathcal{W}^{*}$. Then, by Theorem \ref{thm 2.4}, we have
$\mathcal{R}(\Theta)\subset$~$\mathcal{R}(\Xi).$ The condition $(ii)$ in the result is proved.\\
To show   $ \lambda \Theta^* \Theta \geq \Xi \ \Xi^*$ $(\lambda > 0)$, we consider  upper frame inequality in \eqref{eq 3.7}:
\begin{align*}
 \mathrm{b}_0\|\Theta f\|^2 & \geq \sum_{j,k,m \in \mathbb{Z}}|\langle f,D_{a_j}T_{bk}E_{c_m}\psi\rangle|^2\\
 & =  \sum_{j,k,m \in \mathbb{Z}}|\langle f, \mathcal{W}^{*}e_{j,k,m} \rangle|^2\\
 & = \|\mathcal{W}f \|^2\ \text{ for all} \ f \in L^2(\mathbb{R}).
\end{align*}

This gives $\mathrm{b}_0 \Theta^* \Theta \geq \mathcal{W}^* \mathcal{W}$. That is, $ \lambda \Theta^* \Theta \geq \Xi \ \Xi^*$ $(\lambda  = \mathrm{b}_0 > 0)$. This proves the condition $(i)$ in the result.

Conversely, assume that both conditions $(i)$ and $(ii)$ given in the theorem hold.\\
 We compute
\begin{align*}
\langle \Xi^{*}f,h \rangle &= \left\langle \Xi^{*}f,\sum_{j,k,m \in \mathbb{Z}}a_{j,k,m}e_{j,k,m} \right\rangle\\
&=\sum_{j,k,m \in \mathbb{Z}} \overline{a_{j,k,m}}\langle f, \Xi e_{j,k,m}\rangle\\
&=\sum_{j,k,m \in \mathbb{Z}} \overline{a_{j,k,m}}\langle
f,D_{a_j}T_{bk}E_{c_m}\psi\rangle\\
&=\sum_{j,k,m \in \mathbb{Z}} \overline{\langle
h,e_{j,k,m}\rangle}\langle
f,D_{a_j}T_{bk}E_{c_m}\psi\rangle\\
&=\sum_{j,k,m \in \mathbb{Z}} \langle e_{j,k,m},h\rangle\langle
f,D_{a_j}T_{bk}E_{c_m}\psi\rangle\\
&= \left\langle\sum_{j,k,m \in \mathbb{Z}} \langle
f,D_{a_j}T_{bk}E_{c_m}\psi \rangle e_{j,k,m},h \right\rangle,
\end{align*}
for  all  $f \in L^2(\mathbb{R})$  and for all $h \in \ell^2(\mathbb{Z}^3)$.\\
This gives
\begin{align}
 \Xi^{*}f= \sum_{j,k,m \in \mathbb{Z}} \langle f,
D_{a_j}T_{bk}E_{c_m}\psi \rangle e_{j,k,m} \ \text{for  all} \ f \in L^2(\mathbb{R}).
\end{align}

Therefore, by using $(3.6)$ and the condition $(i)$, we have
\begin{align}\label{eq 3.10}
 \sum_{j,k,m \in \mathbb{Z}} |\langle f,D_{a_j}T_{bk}E_{c_m}\psi\rangle|^2  =\|\Xi^{*}f\|^2 \leq \lambda \|\Theta f\|^2 \ \text{for all} \ f \in L^2(\mathbb{R}) \ (\lambda > 0).
\end{align}

By hypothesis $\mathcal{R}(\Theta)\subset \mathcal{R}(\Xi)$ (see condition $(ii)$). So, by Theorem \ref{thm 2.4}, we can find a positive constant $\beta$  such  that
$\Theta\Theta^{*}\leq \beta \ \Xi \ \Xi^{*}$ (note that $\beta$ is positive, since otherwise $\Theta = O$). Again  by using the
condition $(ii)$,
we have
\begin{align}\label{eq 3.11}
\frac{1}{\beta}  \|\Theta^{*}f\|^2  \leq\|\Xi^{*}f\|^2  & =  \sum_{j,k,m \in \mathbb{Z}} |\langle \Xi^{*}f, e_{j,k,m}\rangle|^2\notag\\
& = \sum_{j,k,m \in \mathbb{Z}} |\langle f, \Xi e_{j,k,m}\rangle|^2\notag\\
&=\sum_{j,k,m \in \mathbb{Z}} |\langle f,D_{a_j}T_{bk}E_{c_m}\psi\rangle|^2 \
\text{for all}\ f \in L^2(\mathbb{R}).
\end{align}
By using \eqref{eq 3.10} and \eqref{eq 3.11}, we conclude that
$\{D_{a_j}T_{bk}E_{c_m}\psi\}_{j,k,m \in \mathbb{Z}} $ is a $\Theta
$-$IWH$ wave packet frame for $L^2(\mathbb{R}).$
\endproof

Djordjevi$\acute{c}$ in \cite{DJ} characterized hyponormal operators  by using the Moore-Penrose inverse of a bounded linear
 operator with a closed range. There may be other  conditions for a bounded linear  operator on a Hilbert space to be hyponormal.
Let $H$ and $K$ be  Hilbert spaces and  $A: H \rightarrow K$ be a bounded linear operator. The Moore-Penrose inverse of $A$ is denoted by $A^{\dag}$, see \cite{BI}. Djordjevi$\acute{c}$ proved the following result by using the Moore-Penrose inverse of a bounded linear operator with a closed range.

\begin{thm}\cite{DJ}
Let $A$ and $AA^* + A^*A$ have closed ranges. Then the following statements are
equivalent:
\begin{enumerate}[$(i)$]
\item $A$ is hyponormal
\item $2AA^*(AA^* +A^*A)^{\dag} AA^* \leq  AA^*$.
\end{enumerate}
\end{thm}
Thus, a bounded linear operator $A$ defined on a Hilbert space is hyponormal if a certain operator  inequality (consisting of adjoint and Moore-Penrose inverse of $A$) is satisfied.  Frame can be used to characterizes  a hyponormal operator on $L^2(\mathbb{R})$. First we define a type of tight frame  (or Parseval frame) in $L^2(\mathbb{R})$. In Definition 3.1, if $\alpha_0 = \beta_0$, then $\{D_{a_j}T_{bk}E_{c_m}\psi\}_{j,k,m \in \mathbb{Z}} $ is not a standard tight  frame, in general. This is the motivation for  new  type of tight frames in $L^2(\mathbb{R})$.
  \begin{defn}
  Let  $\Theta \ne I$ (where $I$ the identity operator on  $L^2(\mathbb{R})$). A  $\Theta$-Hilbert frame   $\{f_n\} \subset \mathcal{H}$ for $\mathcal{H}$ with frame bounds $\alpha_0 = \beta_0$   is  called a \emph{$(\Theta, \alpha_0)$-Hilbert tight  frame}.
\end{defn}

The following theorem characterizes  hyponormal operators  on $L^2(\mathbb{R})$ in terms of $(\Theta, \alpha_0)$- Hilbert tight  frames for $L^2(\mathbb{R})$.
\begin{thm}
A  bounded linear operator $\Theta$ on $L^2(\mathbb{R})$ is hyponormal if and only if there exists  a
$(\Theta, 1)$-Hilbert tight frame for
$L^2(\mathbb{R})$.
\end{thm}
\proof

Assume first that $\Theta$ is a hyponormal operator on $L^2(\mathbb{R})$. Let $\{D_{a_j}T_{bk}E_{c_m}\psi\}_{j,k,m \in \mathbb{Z}}$ be a   tight $I W H$ wave packet frame for $L^2(\mathbb{R})$.\\
Then
\begin{align}\label{eq 3.9}
\sum_{j,k,m \in \mathbb{Z}}|\langle
f,D_{a_j}T_{bk}E_{c_m}\psi\rangle|^2 = \|f\|^2 \  \text{for   all} \
f\in L^2(\mathbb{R}).
\end{align}
Choose $f_n (n \in \mathbb{N}) \leftrightarrow \varphi_{j, k, m} = \Theta(D_{a_j}T_{bk}E_{c_m}\psi), j,k,m \in \mathbb{Z}$.\\
Then, by using  \eqref{eq 3.9} and hyponormality of $\Theta$, we compute
\begin{align}\label{eq 3.5}
\sum_{j,k,m \in \mathbb{Z}}|\langle f, \varphi_{j, k, m}\rangle|^2
&=\sum_{j,k,m \in \mathbb{Z}}|\langle f,\Theta(D_{a_j}T_{bk}E_{c_m}\psi)\rangle|^2\notag\\
 &=\sum_{j,k,m \in \mathbb{Z}}|\langle \Theta^{*}
f,D_{a_j}T_{bk}E_{c_m}\psi\rangle|^2\notag\\
&= \|\Theta^{*} f\|^2\notag\\
&\leq \| \Theta f\|^2 \  \text{for   all} \ f\in
L^2(\mathbb{R}).
\end{align}
For  the lower frame inequality, we compute
\begin{align}\label{eq 3.6}
\|\Theta^{*} f\|^2 & = \sum_{j,k,m \in \mathbb{Z}}|\langle
\Theta^{*} f,D_{a_j}T_{bk}E_{c_m}\psi\rangle|^2\notag\\
&=\sum_{j,k,m \in \mathbb{Z}}|\langle
f,\Theta(D_{a_j}T_{bk}E_{c_m}\psi)\rangle|^2 \notag\\
& = \sum_{j,k,m \in \mathbb{Z}}|\langle f, \varphi_{j, k, m}\rangle|^2\  \text{for   all} \ f\in
L^2(\mathbb{R}).
\end{align}
By using  \eqref{eq 3.5} and \eqref{eq 3.6} we have
\begin{align*}
\|\Theta^{*} f\|^2\leq\sum_{j,k,m \in \mathbb{Z}}|\langle
f, \varphi_{j, k, m}\rangle|^2\leq \|\Theta f\|^2 \
\text{for   all} \ f\in L^2(\mathbb{R}).
\end{align*}
Hence $\{\varphi_{j, k, m}\}_{j,k,m \in \mathbb{Z}}$ is a $(\Theta, 1)$-Hilbert tight  frame for $L^2(\mathbb{R})$.

For the reverse part, suppose that  $\{f_{ n}\}$ is a
$(\Theta, 1)$-Hilbert tight frame for $L^2(\mathbb{R})$.\\
 Then
\begin{align*}
\|\Theta^{*} f\|^2\leq\sum_{n = 1}^{\infty}|\langle
f, f_n\rangle|^2\leq \|\Theta f\|^2 \
\text{for   all} \ f\in L^2(\mathbb{R}).
\end{align*}
This gives $\|\Theta^*f\| \leq \|\Theta f\|$ for all $f \in \mathcal{H}$. Hence $\Theta$  is a hyponormal operator on  $L^2(\mathbb{R})$.
\endproof
%
%

Favier and Zalik proved in \cite{FZ} that the image of a Hilbert frame for $\mathcal{H}$ under a linear homeomorphism is a Hilbert frame for $\mathcal{H}$. They established relation between optimal bounds of a given Hilbert frame and its image (as frame). This is not true for $\Theta$-$IWH$ wave packet frame (see Example \ref{ex 3.13}), in general.  The problem (regarding invariance behaviour as a frame under linear homeomorphism) for $\Theta$-$IWH$ wave packet frames can be solved, provided the given linear homeomorphism commutes with  $\Theta^*$. This is proved in the  following theorem.
\begin{thm}
Let $\mathcal{F} \equiv \{D_{a_j}T_{bk}E_{c_m}\psi\}_{j,k,m \in \mathbb{Z}}$ be  a $\Theta$-$IWH$ wave packet frame for $L^2(\mathbb{R})$ and  $U$   be  a linear homeomorphism on $L^2(\mathbb{R})$ such that  $U$ commutes with  $\Theta^*$. Then, $\mathcal{F}_{U} \equiv \{U(D_{a_j}T_{bk}E_{c_m}\psi)\}_{j,k,m \in \mathbb{Z}}$ is a $\Theta$-$ I W H$ wave packet frame for $L^2(\mathbb{R})$. Furthermore, if  $A_1$ and $B_1$ are optimal bounds of the frame $\mathcal{F}$ and  the pair $(\Theta, U^*)$ is  relatively hyponormal, then the optimal  bounds  $A_2$ and $B_2$ of the frame $\mathcal{F}_{U}$ satisfy the inequalities
\begin{align}\label{eq 3.14}
A_1\|U\|^{-2}\leq A_2\leq A_1\|U^{-1}\|^2  \ ; \ \gamma B_1\|\Theta\|^{-2}\leq
B_2\leq B_1\|U\|^2 \ (\gamma > 0).
\end{align}
\end{thm}
\proof
We compute
\begin{align}\label{eq 3.15}
\sum_{j,k,m \in \mathbb{Z}}|\langle f,U(D_{a_j}T_{bk}E_{c_m}\psi)\rangle|^2
&=\sum_{j,k,m \in \mathbb{Z}}|\langle
U^{*}f,D_{a_j}T_{bk}E_{c_m}\psi\rangle|^2 \notag\\
&\leq B_1\|\Theta U^{*}f\|^2 \notag\\
& = B_1\|U^{*} \Theta f\|^2 \notag\\
&\leq B_1\|U^{*}\|^2\|\Theta f\|^2 \ \text{for all}\ f \in L^2(\mathbb{R}).
\end{align}
By using the fact that $A_1$ is one of the choice for lower $\Theta$-$IWH$ wave packet frame
bound for $\{D_{a_j}T_{bk}E_{c_m}\psi\}_{j,k,m \in \mathbb{Z}}$ and $U$ commutes with $\Theta^*$, we compute
\begin{align}\label{eq 3.16}
\|\Theta^{*}f\|^2 &=\|\Theta^{*}(UU^{-1})f\|^2 \notag\\
&=\|U\Theta^{*}(U^{-1}f)\|^2 \notag\\
&\leq\|U\|^2\|\Theta^{*}(U^{-1}f)\|^2 \notag\\
&\leq\frac{\|U\|^2}{A_1}\sum_{j,k,m \in \mathbb{Z}}|\langle
U^{-1}f,D_{a_j}T_{bk}E_{c_m}\psi\rangle|^2  \notag\\
&=\frac{\|U\|^2}{A_1}\sum_{j,k,m \in \mathbb{Z}}|\langle
UU^{-1}f,U(D_{a_j}T_{bk}E_{c_m}\psi)\rangle|^2\notag\\
&=\frac{\|U\|^2}{A_1}\sum_{j,k,m \in \mathbb{Z}}|\langle
f,U(D_{a_j}T_{bk}E_{c_m}\psi)\rangle|^2.
\end{align}
By using \eqref{eq 3.15} and \eqref{eq 3.16}, we obtain
\begin{align*}
&A_1\|U\|^{-2}\|\Theta^{*}f\|^2\leq\sum_{j,k,m \in \mathbb{Z}}|\langle f,U(D_{a_j}T_{bk}E_{c_m}\psi)\rangle|^2\leq B_1
\|U^*\|^2\|\Theta f\|^2 \ \text{\ for all} \ f\in L^2(\mathbb{R}).
\end{align*}
Hence $\{U(D_{a_j}T_{bk}E_{c_m}\psi)\}_{j,k,m \in \mathbb{Z}}$ is a
$\Theta$-$ I W H$ wave packet frame for
$L^2(\mathbb{R})$ with one of the choice of  frame bounds $A_1\|U\|^{-2}, \
B_1\|U\|^2$.\\
 Since $A_2$ and $B_2$ are best  frame bounds
for $\{U(D_{a_j}T_{bk}E_{c_m}\psi)\}_{j,k,m \in \mathbb{Z}}$, we have
\begin{align}\label{eq 3.17}
 A_1\|U\|^{-2}\leq A_2, \ \  B_2\leq B_1\|U\|^2.
\end{align}
Again  $\{U(D_{a_j}T_{bk}E_{c_m}\psi)\}_{j,k,m \in \mathbb{Z}}$ is a
$\Theta$-$ I W H$ wave packet frame for
$L^2(\mathbb{R})$ with $A_2, B_2$ as one of the choice of  frame bounds.
So, for all $f \in L^2(\mathbb{R}),$ we have
\begin{align}\label{eq 3.18}
A_2\|\Theta^{*}f\|^2\leq\sum_{j,k,m \in \mathbb{Z}}|\langle
f,U(D_{a_j}T_{bk}E_{c_m}\psi)\rangle|^2\leq B_2\|\Theta f\|^2.
\end{align}
For all $f \in L^2(\mathbb{R})$, we have
\begin{align}\label{eq 3.19}
\|\Theta^{*}f\|^2  = \|U^{-1}U\Theta^{*}f\|^2 = \|U^{-1}\Theta^{*}Uf\|^2  \leq\|U^{-1}\|^2\|\Theta^{*}Uf\|^2.
\end{align}
By using \eqref{eq 3.18}, \eqref{eq 3.19} and relative hyponormality of the pair $(\Theta, U^*)$ , we have
\begin{align}\label{eq 3.20}
A_2\|U^{-1}\|^{-2}\|\Theta^{*}f\|^2  &\leq A_2\|\Theta^{*}Uf\|^2 \notag \\
&\leq \sum_{j,k,m \in \mathbb{Z}}|\langle
Uf,U(D_{a_j}T_{bk}E_{c_m}\psi)\rangle|^2 \ \left(=\sum_{j,k,m \in \mathbb{Z}}|\langle f,D_{a_j}T_{bk}E_{c_m}\psi \rangle|^2 \right) \notag\\
&\leq B_2\|\Theta Uf\|^2 \notag\\
&\leq B_2 \|\Theta\|^2 \|Uf\|^2 \notag\\
&\leq  \lambda B_2\|\Theta\|^2\|\Theta f\|^2 \ \text{ for  all} \ f \in
L^2(\mathbb{R}),
\end{align}
where $\lambda$ is a positive constant which appears in the relative hyponormality of the pair $(\Theta, U^*)$.\\
Since  $A_1$ and $B_1$ are the best $\Theta$-$IWH$ wave packet frame bounds for $\{D_{a_j}T_{bk}E_{c_m}\psi\}_{j,k,m \in \mathbb{Z}}$,  by using \eqref{eq 3.20}, we have
\begin{align}\label{eq 3.21}
A_2\|U^{-1}\|^{-2}\leq A_1 , \ \  B_1\leq  \lambda B_2\|\Theta\|^2.
\end{align}
 The  inequalities in \eqref{eq 3.14} are obtained from  \eqref{eq 3.17} and \eqref{eq 3.21}. The result is proved.
\endproof
\begin{rem}
  The condition that the linear homeomorphism $U$ commutes  with  $\Theta^*$ in Theorem 3.10 cannot be relaxed. This is justified  in the following example.
\end{rem}

\begin{exa}\label{ex 3.13}
Consider the multiplication operator  $\Theta:L^2(\mathbb{R})\longrightarrow L^2(\mathbb{R})$  given by
 \begin{align*}
 \Theta(f)=f . \chi_{[0,1]},  \  f \in L^2(\mathbb{R}).
  \end{align*}
Then, $\Theta$ is a bounded linear self-adjoint operator on $L^2(\mathbb{R})$.\\
Choose $a_j=1,c_m=m$ for all $j,m \in \mathbb{Z}$, $b=0$ and $\psi=\chi_{[0,1)}$. Then, $\{D_{a_j}T_{bk}E_{c_m}\psi\}_{j,k,m \in \mathbb{Z}}$ is a
$\Theta$-$ I W H$ wave packet frame for
$L^2(\mathbb{R})$. Indeed, for all $f \in L^2(\mathbb{R})$, we have
\begin{align*}
\sum_{j,k,m \in \mathbb{Z}}|\langle f,D_{a_j}T_{bk}E_{c_m}\psi\rangle|^2
& =\sum_{m \in \mathbb{Z}}|\langle f,E_m\chi_{[0,1)}\rangle|^2\\
& =\sum_{m \in \mathbb{Z}}|\langle f,\Theta(E_m\chi_{[0,1)})\rangle|^2\\
& =\sum_{m \in \mathbb{Z}}|\langle \Theta^* f, E_m\chi_{[0,1)}\rangle|^2\\
& =\|\Theta^* f\|^2\\
& =\|\Theta f\|^2
\end{align*}
Hence $\{D_{a_j}T_{bk}E_{c_m}\psi\}_{j,k,m \in \mathbb{Z}}$ is a $\Theta$-$ I W H$ wave packet frame for
$L^2(\mathbb{R})$. \\

  Choose $U_\circ = T_1 $, the translation operator on   $L^2(\mathbb{R})$, i.e., $U_\circ f(\bullet)= f(\bullet -1) $. Then,
   $ U_\circ$ is a linear homeomorphism on  $L^2(\mathbb{R})$. First we show that  the operator $U_\circ$ and $\Theta^* (=\Theta)$ does not commutes. For this,  we compute
\begin{align}
\Theta^*U_\circ(f)(\gamma)
& =U_\circ(f)(\gamma).\chi_{[0,1)}(\gamma) \notag\\
& =f(\gamma-1).\chi_{[0,1)}(\gamma),\\
\intertext{and}
U_\circ\Theta^*(f)(\gamma)
& =U_\circ(f.\chi_{[0,1)})(\gamma)\notag\\
& =f(\gamma-1).\chi_{[0,1)}(\gamma-1)\notag\\
& =f(\gamma-1).\chi_{[1,2)}(\gamma).
\end{align}
By using $(3.20)$ and $(3.21)$, we conclude that the operators $U_\circ$ and $\Theta^*$ does not commutes.

   Next, we show that the system $ \mathcal{F}_{U_\circ} \equiv  \{U_\circ(D_{a_j}T_{bk}E_{c_m}\psi)\}_{j,k,m \in \mathbb{Z}}$ is not a \break
$\Theta$-$IWH$ wave packet frame for
$L^2(\mathbb{R})$. Let $a_\circ$ and $b_\circ$ be a choice of frame bounds for $ \mathcal{F}_{U_\circ}$.\\
 Then
\begin{align}\label{eq 3.2}
a_o\|\Theta^*f\|^2\leq \sum_{k=1}^{\infty} |\langle f, U_\circ(D_{a_j}T_{bk}E_{c_m}\psi)\rangle|^2 \leq b_o\|\Theta f\|^2 \ \text{for all} \ f \in \mathcal{H}.
\end{align}
Choose $f_\circ = \chi_{[0,1[} \in L^2(\mathbb{R})$. Then, $\|\Theta^* f_\circ \|= 1$.\\
Then, by using lower inequality in  $(3.22)$, we compute
\begin{align*}
a_\circ = a_\circ\|\Theta^*f_\circ\|^2\leq  \sum_{j,k,m \in \mathbb{Z}}|\langle f_\circ, U_\circ(D_{a_j}T_{bk}E_{c_m}\psi)\rangle|^2
 & = \sum_{j,k,m \in \mathbb{Z}}|\langle U_\circ^*f_\circ, D_{a_j}T_{bk}E_{c_m}\psi \rangle|^2\\
 & = \sum_{m \in \mathbb{Z}}|\langle U_\circ^*f_\circ, E_{m}\psi \rangle|^2\\
 & = \|\Theta (U_\circ^* f_\circ)\|^2\\
 & = \|\Theta (\chi_{[-1,0)})\|^2\\
 & =\|\chi_{[-1,0)}.\chi_{[0,1)}\|^2\\
 & =0,
  \end{align*}
a contradiction. Hence  $ \mathcal{F}_{U_\circ}$ is not a $\Theta$-$IWH$ wave packet frame for $L^2(\mathbb{R})$.

\end{exa}

\section{Linear Combinations of $\Theta$-$IWH$ Wave Packet Frames}
Linear combination of frames (or redundant building blocks) is important in applied mathematics.
 Aldroubi in \cite{A} considered the following problem: given a Hilbert frame $\{f_k\}$ for $\mathcal{H}$, define a set of functions $\Phi_{j}$ by taking  linear combinations of the frame elements $f_k$. What are the conditions on the coefficients in the linear combinations, so that the new system $\{\Phi_{j}\}$ constitutes a frame for $\mathcal{H}~?$ More precisely, Aldroubi considered  a linear combination of the
  form
  \begin{align*}\label{eq 4.1}
  \Phi_j = \sum_{k=1}^{\infty}\alpha_{j,k}f_k,   \ (j \in~\mathbb{N})
  \end{align*}
  where $\alpha_{j,k}$ are scalars. Aldroubi  proved sufficient conditions on $\{\alpha_{j,k}\}$  such that
   $\{\Phi_j\}$ constitutes a frame for $\mathcal{H}$. Christensen in  \cite{OC1} gave sufficient conditions which are different from those proved by Aldroubi.
   In this section, we extend some results by Kaushik et al. in \cite{KSV} to  $\Theta$-$IWH$ wave packet frames  for $L^2(\mathbb{R})$.

Let $\{D_{a_j}T_{bk}E_{c_m}\psi\}_{j,k,m \in \mathbb{Z}}$  be a
$\Theta$-$IWH$ wave packet frame for $L^2(\mathbb{R})$. First we  consider a linear combination of the form:
\begin{align}
\Phi_{r,s,t}=\sum_{(j,k,m)\in
\mathbb{I}_{r,s, t}} \alpha_{j,k,m} D_{a_j}T_{bk}E_{c_m}\psi , \
(r,s, t \in \mathbb{Z}),
\end{align}
 where  $\bigcup\limits_{r,s,t \in
\mathbb{Z}}\mathbb{I}_{r,s, t}=\mathbb{Z}\times \mathbb{Z}\times \mathbb{Z}$, \quad $\mathbb{I}_{r,s,t}\bigcap \mathbb{I}_{r',s',t'}= \emptyset$, $(r, s,t)\neq (r', s',t')$  for all $r,s,t,r',s',t' \in~\mathbb{Z}$  and $\alpha_{j,k,m}$ are  scalars. The system
$\{\Phi_{r,s,t}\}_{r,s,t \in \mathbb{Z}}$ is not a $\Theta$-$IWH$ wave packet frame for
$L^2(\mathbb{R})$, in general. This type of combinations under the WH-packet for Gabor system were studied by Kaushik et al.  \cite{KSV}. The following theorem gives  necessary and sufficient conditions for
the system $\{\Phi_{r,s,t}\}_{r,s,t\in \mathbb{Z}}$  to be a $\Theta$-$IWH$ wave packet frame
for $L^2(\mathbb{R})$. This is an adaption of \cite[Theorem 3.5]{KSV}.
\begin{thm}
Let $\Theta$ be a bounded linear operator  on $L^2(\mathbb{R})$ such that $\Theta^*$ is  hyponormal. Assume that  $\{D_{a_j}T_{bk}E_{c_m}\psi\}_{j,k,m \in \mathbb{Z}}$ is a $\Theta$-$IWH$ wave packet frame for $L^2(\mathbb{R})$ and
$\{\Phi_{r,s,t}\}_{r,s,t\in \mathbb{Z}} \subset L^2(\mathbb{R})$ be
the  sequence defined in  \eqref{eq 4.1}. Let \break
$T:\ell^2(\mathbb{Z}^3) \rightarrow \ell^2(\mathbb{Z}^3)$ be a
bounded linear operator such that
\begin{align*}
 T(\{\langle
D_{a_j}T_{bk}E_{c_m}\psi,f\rangle \}_{j,k,m \in \mathbb{Z}})=\{\langle \Phi_{r,s,t},f \rangle \}_{r,s,t\in \mathbb{Z}}, \  f \in L^2(\mathbb{R}).
\end{align*}
 Then, $\{\Phi_{r,s,t}\}_{r,s,t\in \mathbb{Z}}$ is
a $\Theta$-$IWH$ wave packet frame for $L^2(\mathbb{R})$ if and only if there exists a constant
$\lambda>0$ such that
\begin{align}\label{eq 4.2}
\sum\limits_{r,s,t\in \mathbb{Z}}|\langle \Phi_{r,s,t},f\rangle|^2\geq
\lambda\sum\limits_{j,k,m \in \mathbb{Z}}|\langle
D_{a_j}T_{bk}E_{c_m}\psi,f\rangle|^2 \ \text {for  all}  \ f\in
L^2(\mathbb{R}).
\end{align}
\end{thm}
\proof Assume first that  $\{\Phi_{r,s,t}\}_{r,s,t\in \mathbb{Z}}$
is a $\Theta$-$IWH$ wave packet frame for $L^2(\mathbb{R})$ with  frame  bounds $A', B'$. Then, for any
$f\in L^2(\mathbb{R})$, we have
\begin{align}\label{eq 4.3}
\sum\limits_{r,s,t\in \mathbb{Z}}|\langle \Phi_{r,s,t}, f\rangle|^2\geq A'\|\Theta^* f\|^2.
\end{align}
If  $B$ is an upper $\Theta$-$IWH$ wave packet frame bound for $\{D_{a_j}T_{bk}E_{c_m}\psi\}_{j,k,m \in \mathbb{Z}}$, then
\begin{align*}
\sum_{j,k,m \in \mathbb{Z}}|\langle f,
D_{a_j}T_{bk}E_{c_m}\psi\rangle|^2 \leq B \|\Theta f\|^2, \ f\in
L^2(\mathbb{R}).
\end{align*}
i.e.
\begin{align}\label{eq 4.4}
\frac{1}{B}\sum_{j,k,m \in \mathbb{Z}}|\langle f,
D_{a_j}T_{bk}E_{c_m}\psi\rangle|^2 \leq\|\Theta f\|^2, \ f\in
L^2(\mathbb{R}).
\end{align}
Choose  $\lambda=\frac{A'}{B} > 0$. Then, by using hyponormality of $\Theta^*$,  \eqref{eq 4.3}  and
\eqref{eq 4.4}, we have

\begin{align*}
\sum\limits_{r,s,t\in \mathbb{Z}}|\langle \Phi_{r,s,t},f\rangle|^2  & \geq A'\|\Theta^* f\|^2\\
& \geq  A'\|\Theta f\|^2\\
& \geq\lambda\sum_{j,k,m \in \mathbb{Z}}|\langle f, D_{a_j}T_{bk}E_{c_m}\psi\rangle|^2\
\text{for all} \ f\in L^2(\mathbb{R}).
\end{align*}
The inequality given in \eqref{eq 4.2} is proved.\\
For the reverse part, since $\{D_{a_j}T_{bk}E_{c_m}\psi\}_{j,k,m \in \mathbb{Z}}$ is a $\Theta$-$IWH$ wave packet frame for $L^2(\mathbb{R})$. There exist positive
constants  $A, B$  such that
\begin{align}\label{eq 4.5}
A \|\Theta^* f\|^2\leq  \sum_{j,k,m \in \mathbb{Z}}|\langle f,
D_{a_j}T_{bk}E_{c_m}\psi\rangle|^2 \leq B \|\Theta f\|^2 \ \text{for
all} \  f \in L^2(\mathbb{R}).
\end{align}
By using \eqref{eq 4.2} and \eqref{eq 4.5}, we have
\begin{align}\label{eq 4.6}
\sum\limits_{r,s,t\in \mathbb{Z}}|\langle \Phi_{r,s,t},f\rangle|^2&\geq \lambda\sum\limits_{j,k,m \in \mathbb{Z}}|\langle D_{a_j}T_{bk}E_{c_m}\psi, f\rangle|^2 \notag\\
&\geq \lambda A\|\Theta^* f\|^2 \ \text{for  all} \ f \in L^2(\mathbb{R}).
\end{align}
We compute
\begin{align}\label{eq 4.7}
\sum\limits_{r,s,t\in \mathbb{Z}}|\langle \Phi_{r,s,t},f\rangle|^2&=\|\{\langle \Phi_{r,s,t},f\rangle\}_{r,s,t\in \mathbb{Z}}\|^2_{\ell^2(\mathbb{Z}^3)} \notag\\
&=\|T(\{\langle D_{a_j}T_{bk}E_{c_m}\psi,f\rangle\}_{j,k,m \in \mathbb{Z}})\|^2_{\ell^2(\mathbb{Z}^3)} \notag\\
&\leq\|T\|^2\sum_{j,k,m \in \mathbb{Z}}|\langle  D_{a_j}T_{bk}E_{c_m}\psi,f\rangle|^2 \notag\\
&\leq\|T\|^2B\|\Theta f\|^2 \ \text{for  all} \ f \in L^2(\mathbb{R}).
\end{align}
By using \eqref{eq 4.6} and \eqref{eq 4.7}, we conclude that $\{\Phi_{r,s,t}\}_{r,s,t\in \mathbb{Z}}$ is a $\Theta$-$IWH$ wave packet frame for $L^2(\mathbb{R})$.
\endproof

\subsection{The case of finite sum:}
We now  consider a linear combination  of the form $\mathcal{F}_p \equiv\left\{\sum\limits_{s=1}^{p}\alpha_s D_{a_j}T_{bk}E_{c_m}\psi_s
\right\}_{j, k, m \in \mathbb{Z}}$, where $\alpha_1, \alpha_2, ..., \alpha_p$  are nonzero scalars, \break $\psi_s \in L^2(\mathbb{R})$ and  $\{D_{a_j}T_{bk}E_{c_m}\psi_s\}_{j, k, m \in \mathbb{Z}}$  is a  $\Theta$-$IWH$ wave packet
frame for $L^2(\mathbb{R})$ for each $s\in\Lambda_p  = \{1,2,3,..,p\} \ $.  The finite sum $\mathcal{F}_p$
 is not a $\Theta$-$IWH$ wave packet frame for $L^2(\mathbb{R})$, in general.  Kaushik, Singh, and  Virender \cite{KSV} showed that  if some scalar multiple of a  series associated with a  Gabor frame  is dominated by the series associated with the finite sum of Gabor frames, then  the finite sum  constitutes  a Gabor frame for  the underlying space and vice-versa, see Theorem 4.2 of \cite{KSV}. The following theorem  extend   this result in the context of  $\Theta$-$IWH$  wave packet frame for
$L^2(\mathbb{R})$.

\begin{thm}
Assume that $\Theta:L^2(\mathbb{R}) \rightarrow L^2(\mathbb{R})$ is a bounded linear operator such that  $\Theta^*$ is  hyponormal.
Let $\{D_{a_j}T_{bk}E_{c_m}\psi_s\}_{{j, k, m \in \mathbb{Z}}\atop{s\in\Lambda_p}}$
 be a finite family of $\Theta$-$IWH$ frames for $L^2(\mathbb{R})$. Then,
$\mathcal{F}_p \equiv \left\{\sum\limits_{s=1}^{p}\alpha_s D_{a_j}T_{bk}E_{c_m}\psi_s
\right\}_{j, k, m \in \mathbb{Z}}$ is a $\Theta$-$IWH$ wave packet frame for $L^2(\mathbb{R})$ if
and only if there exists $\mu>0$ and some $\xi \in \Lambda_p$ such
that
\begin{align*}
\mu\sum_{j, k, m \in \mathbb{Z}}|\langle
D_{a_j}T_{bk}E_{c_m}\psi_\xi,f\rangle|^2\leq\sum_{j, k, m \in
\mathbb{Z}}\left|\left\langle \  \sum\limits_{s=1}^p \alpha_s
D_{a_j}T_{bk}E_{c_m}\psi_s,f \right\rangle \right|^2, \  f\in L^2(\mathbb{R})
\end{align*}
for any finite  sequence of scalars $\{\alpha_s\}$.
\end{thm}
\proof Let $A_\xi, B_\xi$ be  frame  bounds for
 $\Theta$-$IWH$ wave packet frame $\{D_{a_j}T_{bk}E_{c_m}\psi_\xi\}_{j, k,  m \in \mathbb{Z}}$ for $L^2(\mathbb{R})$ $(1\leq \xi \leq p)$.\\
 Then
\begin{align}
\mu A_\xi \|\Theta^* f\|^2&\leq\mu\sum_{j, k, m \in \mathbb{Z}}|\langle  D_{a_j}T_{bk}E_{c_m}\psi_\xi, f\rangle|^2 \notag\\
&\leq\sum_{j, k, m \in \mathbb{Z}}\left| \left \langle  \sum\limits_{s=1}^p
\alpha_s D_{a_j}T_{bk}E_{c_m}\psi_s, f \right \rangle \right|^2, \ f\in L^2(\mathbb{R}).
\end{align}
Thus, the lower frame condition is satisfied for the finite system $\mathcal{F}_p$.\\

For the upper  frame condition, we compute
\begin{align}
\sum_{j, k, m \in \mathbb{Z}}\left|\left\langle \ \sum\limits_{s=1}^p \alpha_s D_{a_j}T_{bk}E_{c_m}\psi_s, f \right\rangle \right|^2
&=\sum_{j, k, m \in \mathbb{Z}}\left|\sum\limits_{s=1}^p \alpha_s\langle   D_{a_j}T_{bk}E_{c_m}\psi_s,f\rangle \right|^2 \notag\\
&\leq\sum_{j, k, m \in \mathbb{Z}}\left[\sum\limits_{s=1}^p \left| \alpha_s\langle  D_{a_j}T_{bk}E_{c_m}\psi_s,f\rangle \right|\right]^2 \notag\\
&\leq \sum\limits_{s=1}^p \left( |\alpha_s|^2\sum_{j, k, m \in \mathbb{Z}}| \langle D_{a_j}T_{bk}E_{c_m}\psi_s,f\rangle|^2 \right)\notag\\
&\leq \left( p \max_{1 \leq s \leq p}|\alpha_s|^2\sum\limits_{s=1}^p B_s \right) \|\Theta f\|^2, \ f\in
L^2(\mathbb{R}).
\end{align}
By $(4.8)$ and $(4.9)$, we conclude that the finite sum $\mathcal{F}_p$
 is a $\Theta$-$IWH$ wave packet frame for $L^2(\mathbb{R})$.\\
Conversely, assume that the finite sum  $\mathcal{F}_p$
 is  a $\Theta$-$IWH$ wave packet frame for $L^2(\mathbb{R})$
with frame bounds $A,B$. Then, for all $f\in L^2(\mathbb{R})$, we have
\begin{align}
A\|\Theta^* f\|^2\leq\sum_{j, k, m \in \mathbb{Z}}\left|\left\langle
\sum\limits_{s=1}^p \alpha_s D_{a_j}T_{bk}E_{c_m}\psi_s,
f \right \rangle \right|^2.
\end{align}
If  $B_\xi$ is an upper frame bound for $\{D_{a_j}T_{bk}E_{c_m}\psi_{\xi}\}_{j, k, m \in \mathbb{Z}}$, then
\begin{align}
\frac{1}{B_\xi}\sum_{j, k, m \in \mathbb{Z}}|\langle f,
D_{a_j}T_{bk}E_{c_m}\psi_{\xi}\rangle|^2 \leq \|\Theta f\|^2, \ f\in
L^2(\mathbb{R}).
\end{align}
Choose $\mu=\frac{A}{B_\xi} > 0$. Then,
using hyponormality of $\Theta^*$, $(4.10)$ and $(4.11)$ we have
\begin{align*}
\mu \sum_{j, k, m \in \mathbb{Z}}|\langle
D_{a_j}T_{bk}E_{c_m}\psi_\xi, f\rangle|^2&\leq A\|\Theta f\|^2\\
&\leq A\|\Theta^* f\|^2\\
&\leq\sum_{j, k, m \in\mathbb{Z}}\left| \left\langle  \sum\limits_{s=1}^p \alpha_s
D_{a_j}T_{bk}E_{c_m}\psi_s, f\right\rangle \right|^2,\ f\in L^2(\mathbb{R}).
\end{align*}
The theorem is proved.
\endproof

\textbf{Application}: The following example gives an application of Theorem 4.2.
\begin{exa}
Let $\Theta: L^2(\mathbb{R}) \rightarrow L^2(\mathbb{R})$ be the modulation operator. That is,  $\Theta f(t) = e^{2\pi ibt}f(t)$, where $b \in \mathbb{R}$ is fixed. Then,  $\Theta^*$ is  hyponormal  on $L^2(\mathbb{R})$.\\
   Choose $\psi=\chi_{[0,1[}$,  $a_j = 1, \ b =1, c_m = m$ for all $j,m\in\mathbb{Z}$  and  $\psi_s=\psi$  for all $s\in\Lambda_p$.
Then,  for any nonzero scalars $\alpha_1, \alpha_2, ..., \alpha_p$ with $\sum_{s=1}^p\alpha_s \ne 0$, we have
\begin{align*}
\sum_{j,k,m \in
\mathbb{Z}}\left| \left\langle  \sum\limits_{s=1}^p \alpha_s
D_{a_j}T_{bk}E_{c_m}\psi_s,f\right\rangle \right|^2 &=\sum_{j,k,m\in\mathbb{Z}}\left|\left\langle\sum_{s=1}^p \alpha_s D_{a_j}T_{bk}E_{c_m} \psi,f \right\rangle \right|^2\\
& = |\sum_{s=1}^p\alpha_s|^2\sum_{j,k,m\in \mathbb{Z}}|\langle D_{a_j}T_{bk}E_{c_m} \psi_\xi, f\rangle|^2, \ \  f\in L^2(\mathbb{R}),
\end{align*}
where $\psi_\xi = \chi_{[0,1)}$.\\
Choose $\mu = |\sum_{s=1}^p\alpha_s|^2 > 0$.\\
Then
\begin{align*}
\mu \sum_{j,k,m \in \mathbb{Z}}|\langle
D_{a_j}T_{bk}E_{c_m}\psi_\xi,f\rangle|^2=\sum_{j,k,m \in
\mathbb{Z}}\left| \left\langle  \sum\limits_{s=1}^p \alpha_s
D_{a_j}T_{bk}E_{c_m}\psi_s,f\right\rangle \right|^2 \ \text{for all} \  f\in L^2(\mathbb{R}).
\end{align*}

By Theorem 4.2, the finite system  $\mathcal{F}_p$ is a $\Theta$-$IWH$ wave packet frame for $L^2(\mathbb{R})$.
\end{exa}

\begin{center}
  This work is jointly with A. K. Sah and Deepshikha
\end{center}

$$\textbf{Acknowledgement}$$
Lalit is   supported by R $\&$ D Doctoral Research Programme, University of Delhi,
Delhi-110007, India. Grant No. : RC/2014/6820.

\bibliographystyle{amsplain}

\mbox{}

\mbox{}

\end{document}